\documentclass{amsart}

\usepackage{amsfonts}
\usepackage{amssymb, latexsym, amsmath, pb-diagram}
\usepackage{pstricks-add}
\usepackage{lamsarrow, pb-lams}
\usepackage{multicol}
\usepackage{bm}
\usepackage[mathscr]{eucal}
\usepackage{xy}
\usepackage{epic,eepic}
\xyoption{all}

\def\w{\widetilde}

\def\nr{\refstepcounter{thm}\thethm}
\def\cf{{\it cf.}\ }

\begin{document}
\title[On $H^*(S(n,k)$)]{The cohomology of $(S(n,k))$ relevant to  Morava stabilizer algebra}
\author[L.~Chen, X.~Wang \& X. Zhao]{Liman Chen, Xiangjun Wang and Xuezhi Zhao}
\date{}
\address{School of Mathematical Science, Nankai University, Tianjin 300071,
P.~R.~China}
\email{chenlimanstar1@163.com}
\address{School of Mathematical Science and LPMC, Nankai University, Tianjin 300071,
P.~R.~China}
\email{xjwang@nankai.edu.cn}
\address{Department of Mathematics \& Institute of mathematics and interdisciplinary
science, Capital Normal University, Beijing 100048, P.~R.~China}
\email{zhaoxve@mail.cnu.edu.cn}
\keywords{Morava n-th
K-theory, Morava stabilizer algebra, May spectral sequence.}
\subjclass[1991]{Primary 55Q99; Secondary 55Q51}
\thanks{The second author was supported by NSFC grant No. 11471167 and SRFDP No.20120031110025}
\maketitle

\begin{abstract}
In this paper we redefine a increasing filtration on the the Hopf algebra $S(n,k)$,
From which we get a spectral sequence called May spectral sequence. As an application
we computed $H^{*,*}S(n,n)$ at prime 2, $H^{*,*}S(3,2)$ at prime 3 and $H^{*,*}S(4,2)$
at prime $p\geqslant 5$ 
\end{abstract}

\section{Introduction}
    In stable homotopy theory, the ``chromatic'' point of view plays an important role
(\cf \cite{Ghmr, RavG, Rav84}).
Fix a prime $p$. Let $E(n)_*$, $n\geqslant 0$ be the Johnson-Wilson
homology theories and let $L_n$ be localization functor with respect
to $E(n)_*$. Then there are natural transformations
$L_nX\longrightarrow L_{n-1}X$, and the chromatic tower
\[
\xymatrix{
\cdots \ar[r] & L_nX \ar[r] & L_{n-1}X \ar[r] & \cdots \ar[r] &
L_2X \ar[r] & L_1X \ar[r] & L_0X.
}
\]
By the Hopkins-Ravenel chromatic convergence theorem, the homotopy
inverse limit of this tower is the
$p$-localization of $X$
\[
X\longrightarrow \mbox{Holim}L_nX.
\]
 Thus the homotopy groups
$\pi_*(L_nX)$ is the part of homotopy groups $\pi_*(X)$ one
could see from $E(n)_*$.

  To determine the homotopy groups $\pi_*(L_nX)$, one has the Adams-Novikov
spectral sequence based on the Brown-Peterson spectrum $BP$, whose $E_2-$term
is
\[
E_2^{s,t}=Ext^{s,t}_{BP_*BP}(BP_*, BP_*(L_nX)).
\]
(\cf \cite{Adams, Mrw, RavG, Rav84})

    To determine the Adams-Novikov $E_2$-term $Ext^{s,t}_{BP_*BP}(BP_*, BP_*(L_nX))$
one has the Bockstein spectral sequence. This is an argument based on
the cohomology of the Morava stabilizer algebra $S(n)$ at each prime $p$ (\cf \cite{RavG, ShW, ShW1, W}).
Here the Hopf algebra $S(n)$ is defined as
$$S(n)=Z/p\otimes_{K(n)_*}K(n)_*K(n)\otimes_{K(n)_*}Z/p,$$
where $K(n)_*=Z/p[v_n, v_n^{-1}]$,
 $$K(n)_*K(n)=K(n)_*\otimes_{BP_*}BP_*BP\otimes_{BP_*}K(n)_*=K(n)_*[t_1, t_2, \cdots]/(v_nt_s^{p^n}-v_n^{p^s}t_s),$$
$K(n)_*$ acts on $Z/p$ by sending $v_n$ to $1$. Thus
$$S(n)=Z/p[t_1, t_2, \cdots, t_s, \cdots]/(t_s^{p^n}-t_s).$$

    We write $S(n,k)=S(n)/(t_j:j<k)=Z/p[t_k, t_{k+1}, \cdots t_s, \cdots]/(t_s^{p^n}-t_s)$.
The Hopf algebra structure of $S(n)$ determines that of $S(n,k)$, while $S(n,1)=S(n)$.
Let $V(n-1)$ and $T(k-1)$  denote the Smith-Toda spectra and the
Ravenel spectra respectively characterized by
\begin{align*}
BP_*V(n-1) & = BP_*/I_{n}=BP_*/(p,v_1,\cdots, v_{n-1}) & & \mbox{and}\\
 BP_*T(k-1) & =BP_*[t_1, t_2,\cdots, t_{k-1}].
\end{align*}
If $L_nV(n-1)\wedge T(k-1)$ exist, (although $V(n-1)$ does not exist (\cf \cite{Lee}),
but $V(n-1)\wedge T(k-1)$ might exist),
then by the change of rings theorem, the $E_2$-term of the
Adams-Novikov spectral sequence converging to $\pi_*(L_nV(n-1)\wedge T(k-1))$ is
\begin{align*}
      & Ext_{BP_*BP}^{s,t}(BP_*, BP_*(L_nV(n-1)\wedge T(k-1)))\\
\cong & Ext_{S(n,k)}^{s,t}(Z/p, Z/p) \otimes K(n)_*[v_{n+1}, \cdots, v_{n+k-1}].
\end{align*}

   In this paper, we will use $H^{s,*}S(n,k)$ to denote the $Ext$ groups
$Ext_{S(n,k)}^{s,*}(Z/p, Z/p)$. In \cite{Henn, Rav77}, Ravenel and Henn determined $H^{s,*}S(1)$,
$H^{s,*}S(2)$ at all primes, and $H^{s,*}S(3)$ at the odd primes $p\geqslant 5$.
$H^{s,*}S(n,k)$ is known from \cite{RavG} for $k\geqslant n$ at odd primes and $k>n$ at the prime 2 .
In \cite{ShT} Shimomura and Tokashiki computed $H^{s,*}S(n,n-1)$ at odd primes $p>3$.
In this paper we will be concentrated on the case $k\leqslant n$.

    Consider the cohomology of the Hopf algebra $S(n,k)$ at all primes.
In section 2 of this paper,
we follow Ravenel's ideal (\cf \cite{RavG} 3.2.5 Theorem), redefined the May filtration
in $S(n,k)$ and its cobar complex $C^{s,*}S(n,k)$.
This filtration induces a spectral sequence so called May spectral
sequence $\{E_r^{s,*,M}(n,k), d_r\}$ that converges to
$H^{s,*}S(n,k)$. Then in section 3 we prove that the $E_2$-term of the
May spectral sequence is
isomorphic to the cohomology of
\[
\left\{E[h_{i,j}|k\leqslant i\leqslant s_0, j\in Z/n]
    \otimes P[b_{i,j}|k\leqslant i\leqslant s_0-n, j\in Z/n], d_1\right\}
\]
where $s_0=max\left\{\left[{\displaystyle \frac{2pn+p-2}{2(p-1)}}\right], n+k-1\right\}$
and $\left[{\displaystyle \frac{2pn+p-2}{2(p-1)}}\right]$ is the integer part of
${\displaystyle \frac{2pn+p-2}{2(p-1)}}$.
In particular, if
\[
n+k-1\geqslant \left[{\displaystyle \frac{2pn+p-2}{2(p-1)}}\right],
\]
the May's $E_2$-term becomes the cohomology of
\[
\left\{E[h_{i,j}|k\leqslant i\leqslant n+k-1,\ j\in Z/n],\ d_1\right\}.
\]
The homological dimension of each element is given by
\begin{align*}
s(h_{i,j})= & 1, & s(b_{i,j})=& 2.
\end{align*}
For the May differentials, one has $d_r:E_r^{s,t,M}S(n,k)\longrightarrow
E_r^{s+1,t,M-r}S(n,k)$ and if $x\in E_r^{s,*,*}S(n,k)$ then
\[
d_r(x\cdot y)=d_r(x)\cdot y+(-1)^{s}x\cdot d_r(y).
\]
The first May differential $d_1$ is given by
\begin{align*}
 d_1(h_{i,j})= & -\sum_{k\leqslant m\leqslant i-k}h_{m,j}h_{i-m, j+m}, &
 \mbox{and} \hspace{10mm} d_1(b_{i,j})&= 0.
\end{align*}
We analyze the higher May differentials and give a collapse
theorem in section 4. As an
consequence we compute the cohomology of $S(n,n)$ at the prime 2,
$S(3,2)$ at the prime 3 and $S(4,2)$ at the prime $p\geqslant 5$ in section 5.

\section{The May spectral sequence}

   Let $p$ be a prime, $BP_*=Z_{(p)}[v_1, v_2, \cdots ]$
and $BP_*BP=BP_*[t_1, t_2, \cdots ]$.
For the Hazewinkel's generators
described inductively by
$v_s=pm_s-{\displaystyle \sum_{i=1}^{s-1}v_{s-i}^{p^i}m_i}$ (\cf \cite{Haz} and \cite{Mrw} 1.2), the diagonal map
$\Delta: BP_*BP\rightarrow BP_*BP\otimes_{BP_*}BP_*BP$ is given by
$$\sum_{i+j=s}m_i(\Delta t_j)=\sum_{i+j+k=s}m_it_j^{p^i}\otimes t_k^{p^{i+j}}.$$
One can easily prove that
\begin{align*}
\Delta(t_1) = & t_1\otimes 1+ 1\otimes t_1 &  & \mbox{and}\\
\Delta(t_2) = & \sum_{i+j=2}t_i\otimes t_j^{p^i} - v_1b_{1,0},
\end{align*}
where $p\cdot b_{1,0}=\Delta(t_1^p)-t_1^p\otimes 1-1\otimes t_1^p$.  Inductively
define
$$p\cdot b_{s,k-1}=\Delta(t_s^{p^k})-\sum_{i+j=s}t_i^{p^k}\otimes t_j^{p^{i+k}}
     +\sum_{0<i<s} v_i^{p^k}b_{s-i, k+i-1},$$
one has
\begin{align}\tag\nr
\Delta(t_{s+1})= \sum_{i+j=s+1}t_i\otimes t_j^{p^i}-\sum_{0<i<s+1}v_ib_{s+1-i, i-1}.
\end{align}

Thus for the $n$-th Morava K-theory $K(n)=Z/p[v_n, v_n^{-1}]$, the Hopf algebra
\[
K(n)_*K(n)=K(n)_*\otimes_{BP_*}BP_*BP\otimes_{BP_*}K(n)_*
\] is isomorphic to
\[
K(n)_*K(n)=Z/p[t_1,\ t_2,\ \cdots,\ t_s,\ \cdots ]/(v_nt_s^{p^n}-v_n^{p^s}t_s).
\]
And $S(n)=Z/p\otimes_{K(n)_*}K(n)_*K(n)\otimes_{K(n)_*}Z/p$ is isomorphic to
\[
S(n)=Z/p[t_1, t_2, \cdots, t_n, t_{n+1}, \cdots ]/(t_s^{p^n}-t_s).
\]
The inner degree of $t_s$ in $S(n)$ is
\begin{align*}
|t_s| & \equiv 2(p-1)(1+p+\cdots +p^{s-1}) & & \mbox{$mod$ $2(p-1)(1+p+\cdots +p^{n-1})$}
\end{align*}
because $v_n$ is sent to $1$.  The structure map $\Delta: S(n)\longrightarrow S(n)\otimes S(n)$
acts on $t_s$ as follows
\begin{align*}
\Delta(t_s) & =\sum_{0\leqslant i\leqslant s}t_i\otimes t_{s-i}^{p^i} & &\mbox{for $s\leqslant n$} \\
\Delta(t_s) & =\sum_{0\leqslant i\leqslant s}t_i\otimes t_{s-i}^{p^i} -b_{s-n, n-1}
                          & & \mbox{for $s>n$}
\end{align*}
where $b_{i,j}=\sum_{0<m<p}{p \choose m}/p\ \cdot\
   t_{i}^{mp^{j}}\otimes t_{i}^{(p-m)p^{j}}$ at odd primes and
$b_{i,j}=t_{i}^{2^{j}}\otimes t_{i}^{2^{j}}$ at the prime 2.
For the integer $k\geqslant 1$, let $S(n,k)=S(n)/(t_s|s<k)$. We have
$$S(n,k)=Z/p[t_k, t_{k+1}, \cdots, t_{n+k}, t_{n+k+1}, \cdots ]/(t_s^{p^n}-t_s),$$
the structure map $\Delta: S(n,k)\longrightarrow S(n,k)\otimes S(n,k)$ acts on
$t_s$ as
\begin{align}\notag
\Delta(t_s) & =1\otimes t_s+\sum_{k\leqslant i\leqslant s-k} t_i\otimes t_{s-i}^{p^i}+t_s\otimes 1 & &\mbox{for $s\leqslant n+k-1$,}
        \notag \\
\Delta(t_s) & =1\otimes t_s+\sum_{k\leqslant i\leqslant s-k} t_i\otimes t_{s-i}^{p^i}+t_s\otimes 1 \tag{\nr}
            -b_{s-n, n-1}
                          & & \mbox{for $s\geqslant n+k$.}
\end{align}

    In the resulting May spectral sequence, we want to have the 0-th May differential is
\begin{align*}
d_0(t_s^{p^j}) & = 0
\end{align*}
and the first May differential is given by
\begin{align*}
d_1(t_s) & = t_k\otimes t_{s-k}^{p^k} +t_{k+1}\otimes t_{s-k-1}^{p^{k+1}}+\cdots + t_{s-k}\otimes t_k^{p^{s-k}} & &
 \mbox{for $s\leqslant n+k-1$,}
\end{align*}
and for $s\geqslant n+k$
\begin{align*}
d_1(t_s) & = \left\{\begin{array}{llll}
                    t_k\otimes t_{s-k}^{p^k} +\cdots + t_{s-k}\otimes t_k^{p^{s-k}}&
                    \mbox{if the May flitration $M(t_k\otimes t_{s-k}^{p^k})>M(b_{s-n,n-1})$,}\\
                    -b_{s-n,n-1} &
                    \mbox{if the May flitration $M(t_k\otimes t_{s-k}^{p^k})\leqslant M(b_{s-n,n-1})$.}
                    \end{array}\right.
\end{align*}
So we define the May filtration on $S(n,k)$ as:

{\bf Definition \nr} {\it In the Hopf algebra $S(n,k)$, we define May filtration $M$ as follows:
\begin{enumerate}
\item For $k\leqslant s\leqslant n+k-1$, set the May filtration of $t_s^{p^j}$ as
        $M(t_s^{p^j})=2s-1$.
\item For $n+k\leqslant s$, inductively set the May filtration of $t_{s}^{p^j}$ as
       $$M(t_{s}^{p^j})=max\{2s-1,\ \ p M(t_{s-n}^{p^{j+n-1}})+1\}.$$
\item For the monomial $t_s^j$, express $j$ by the $p-$adic number as $j=j_0+j_1p+\cdots +j_mp^m$,
    where $0\leqslant j_i<p$. Set the May filtration of $t_s^j$ as
      $$M(t_s^j)=\sum_{0\leqslant i\leqslant m}j_iM(t_s^{p^i}).$$
\item And for $t_{s_1}^{j_1}\cdot t_{s_2}^{j_2}\cdots t_{s_m}^{j_m}$, where $s_i\not=s_j$
    define its May filtration as
      $$M(t_{s_1}^{j_1}\cdot t_{s_2}^{j_2}\cdots t_{s_m}^{j_m})=\sum_{1\leqslant i\leqslant m}M(t_{s_i}^{j_i}).$$
\end{enumerate}}

{\bf Lemma \nr} {\it Let $s_0=max\left\{\left[\displaystyle{\frac{2pn+p-2}{2(p-1)}}\right],\ n+k-1\right\}$
where $\left[\displaystyle{\frac{2pn+p-2}{2(p-1)}}\right]$ is the integer part of $\displaystyle{\frac{2pn+p-2}{2(p-1)}}$.
Then the May filtration of $t_s^{p^j}$ satisfies
\begin{enumerate}
\item $M(t_{s}^{p^j})>M(t_{s-1}^{p^j})+1$ and
\item For $s\leqslant s_0$, the May filtration $M(t_s^{p^j})=2s-1$.
\item For $s>s_0$, $pM(t_{s-n}^{p^j})+1\geqslant 2s-1$ and the May filtration $M(t_s^{p^j})= pM(t_{s-n}^{p^{j}})+1$.
\end{enumerate}}

\begin{proof}
    1) If $s_0=max\left\{\left[\frac{2pn+p-2}{2(p-1)}\right],\ n+k-1\right\}=n+k-1$.
From its definition, we see that for $s\leqslant n+k-1=s_0$, the May filtration
of $t_s^{p^j}$ is $2s-1$ and $M(t_s^{p^j})>M(t_{s-1}^{p^j})+1$.
From $n+k-1\geqslant \left[\frac{2pn+p-2}{2(p-1)}\right]$, one sees that
\begin{align*}
s_0+1=n+k > & \frac{2pn+p-2}{2(p-1)} & & \mbox{and} & p(2k-1)+1> & 2(n+k)-1.
\end{align*}
Thus from $M(t_k^{p^j})=2k-1$, one knows that the May filtration of $t_{n+k}^{p^j}$
is $pM(t_k^{p^j})+1$ and
\[
M(t_{n+k}^{p^j})=p(2k-1)+1>2(n+k)-1=M(t_{n+k-1}^{p^j})+1.
\]
Inductively suppose that $M(t_s^{p^j})>  M(t_{s-1}^{p^j})+1$ and for $s_0<s\leqslant m$,
\[
pM(t_{s-n}^{p^j})+1>2s-1,
\]
so the May filtration $M(t_s^{p^j})=  pM(t_{s-n}^{p^j})+1$. Then
from $M(t_{m+1-n}^{p^j})>M(t_{m-n}^{p^j})+1$ one get
\begin{align*}
pM(t_{m+1-n}^{p^j})+1 > & p\left(M(t_{m-n}^{p^j})+1\right) +1=pM(t_{m-n}^{p^j})+p+1\\
                      > & 2m-1 +p \geqslant 2(m+1)-1.
\end{align*}
The May filtration of $t_{m+1}^{p^j}$ is $pM(t_{m+1-n}^{p^j})+1$.

If $s_0=\left[\displaystyle{\frac{2pn+p-2}{2(p-1)}}\right]> n+k-1$, then for $k\leqslant s\leqslant s_0$,
$s\leqslant \displaystyle{\frac{2pn+p-2}{2(p-1)}}$. This implies
\[
p\left(2(s-n)-1\right)+1\leqslant 2s-1.
\]
From $\displaystyle{\frac{2pn+p-2}{2(p-1)}}\leqslant 2n$ we see that
$s-n\leqslant n\leqslant n+k-1$. Thus the May filtration $M(t_{s-n}^{p^j})=2(s-n)-1$
and $pM(t_{s-n}^{p^j})+1<2s-1$.
This implies that the May filtration
of $t_s^{p^j}$ is $2s-1$ and $M(t_s^{p^j})>M(t_{s-1}^{p^j})+1$.

   Notice that $s_0+1=\displaystyle{\left[\frac{2pn+p-2}{2(p-1)}\right]}+1>\displaystyle{\frac{2pn+p-2}{2(p-1)}}$,
this implies
\[
p(2(s_0+1-n)-1)+1>2(s_0+1)-1.
\]
The May filtration of $t_{s_0+1-n}^{p^j}$ is $2(s_0+1-n)-1$, so the May filtration
\[
M(t_{s_0+1}^{p^j})=pM(t_{s_0+1-n}^{p^j})+1.
\]
Similarly, by induction we get the Lemma.
\end{proof}

{\bf Example:} The May filtration in $S(4,2)$ is given by:
\[
\begin{array}{cccccccccccc}
    & t_2 & t_3 & t_4 & t_5 & t_6 & t_7 & t_8 & t_9 & t_{10} & \cdots \\
p=2 & 3   & 5   & 7   & 9   & 11  & 13  & 15  & 19  & 23     & \cdots \\
p=3 & 3   & 5   & 7   & 9   & 11  & 16  & 22  & 28  & 34     & \cdots \\
p\geqslant 5 & 3 & 5 & 7 & 9 & 3p+1 & 5p+1 & 7p+1 & 9p+1 & p(3p+1)+1 & \cdots
\end{array}
\]

   Let $F^{*,M}(n,k)$ be the sub-module of $S(n,k)$ generated by the elements with May
filtration $\leqslant M$. Set
$E^{*,M}(n,k)=F^{*,M}(n,k)/F^{*,M-1}(n,k)$.
One can see from Lemma 2.4 that
\begin{align}\tag\nr
E^{*,*}(n,k) \cong \bigotimes_{k\leqslant s}T[t_s^{p^j}|j\in Z/n]
\end{align}
is a bigraded Hopf algebra, where $T[\ \ \ ]$ denote the truncated polynomial algebra of height
$p$ on the indicated generators. The structure map
\[
\Delta: E^{*,*}(n,k) \longrightarrow E^{*,*}(n,k) \otimes E^{*,*}(n,k)
\]
acts the the generators $t_s^{p^j}$ as $\Delta(t_s^{p^j})=1\otimes t_s^{p^j}+t_s^{p^j}\otimes 1$.

    Let $C^{s,t}S(n,k)=\otimes^s\bar{S}(n,k)$ denote the cobar construction of $S(n,k)$
where $\bar{S}(n,k)=Ker\ \epsilon$ denote the augmentation ideal of $S(n,k)$. The differential
\[
d: C^{s,t}S(n,k)\longrightarrow C^{s+1,t}S(n,k)
\]
is given on the generators as
\begin{align}\notag
  & d(\alpha_1\otimes \cdots \otimes \alpha_s)\\
= & \sum_{1\leqslant i\leqslant s}(-1)^{i}\alpha_1\otimes \cdots
\otimes (\Delta(\alpha_i)-\alpha_i\otimes 1-1\otimes \alpha_i) \otimes \cdots
\otimes \alpha_s. \tag{\nr}
\end{align}
In general, the generator $\alpha_1\otimes \alpha_2 \otimes \cdots \otimes \alpha_s$
of $C^{s,t}S(n,k)$ is denoted by $[\alpha_1|\alpha_2|\cdots|\alpha_s]$.
For the generator $[\alpha_1|\alpha_2|\cdots|\alpha_s]$, define its May
filtration as
$$M([\alpha_1|\alpha_2|\cdots|\alpha_s])=M(\alpha_1)+M(\alpha_2)+\cdots +M(\alpha_s).$$

    Let $FC^{*,*,M}S(n,k)$ denote the sub-complex of $C^{*,*}S(n,k)$ generated by the elements
with May filtration $\leqslant M$.  Then we get a short exact sequence
\begin{align}\tag{\nr}
0\longrightarrow FC^{*,*,M-1}S(n,k) \longrightarrow FC^{*,*,M}S(n,k)
  & \longrightarrow E_0^{*,*,M}S(n,k) \longrightarrow 0
\end{align}
of cochain complexes. The cochain complex
$$E_0^{*,*,M}S(n,k)=FC^{*,*,M}S(n,k)/FC^{*,*,M-1}S(n,k)$$
is isomorphic to the cobar complex of $E^{*,*}(n,k)$ given in (2.5).
Let $E_1^{*,*,M}S(n,k)$ be the homology of $(E_0^{*, *, M}S(n,k), d_0)$. Then
(2.7) gives rise to a spectral sequence (so called the May spectral sequence)
\[
\{E_r^{s,t,M}(n,k), d_r\}
\]
that converges to
\[
H^{s,t}(C^{*,t}S(n,k), d)= Ext_{S(n,k)}^{s,t}(Z/p, Z/p).
\]

{\bf Theorem \nr} {\it For $k\leqslant n$ the Hopf algebra $S(n,k)$ can be given an increasing filtration
as in Definition (2.3). The associated bigraded Hopf algebra $E^{*,M}(n,k)$ is primitively generated with the algebra
structure of (2.5). In the associated spectral sequence,
the $E_1$-term $E_1^{s,t,M}S(n,k)$ is isomorphic to
$$E[h_{i,j}|k\leqslant i,j\in Z/n]\otimes P[b_{i,j}|k\leqslant i, j\in Z/n].$$
The homological dimension of each element is given by
$s(h_{i,j})=1$, $s(b_{i,j})=2$ and the degree is given by
\begin{align*}
h_{i,j}\in & E_1^{1,2(p^i-1)p^j,*}(n,k),\\
b_{i,j}\in & E_1^{2,2(p^i-1)p^{j+1},*}(n,k)
\end{align*}
here $h_{i,j}$ corresponds to $t_i^{p^j}$ and
$b_{i,j}$ corresponds to
$\sum {p \choose m}/pt_i^{mp^j} \otimes t_i^{(p-m)p^j}$.
One has $d_r:E_r^{s,t,M}S(n,k)\longrightarrow
E_r^{s+1,t,M-r}S(n,k)$ and if $x\in E_r^{s,*,*}S(n,k)$ then
$$d_r(x\cdot y)=d_r(x)\cdot y+(-1)^{s}x\cdot d_r(y).$$
In the $E_1$-term of this spectral sequence, we have the
following relations:
\begin{align*}
h_{i,j}\cdot h_{i_1,j_1} & =-h_{i_1,j_1}h_{i,j} & h_{i,j}\cdot b_{i_1,j_1} & =b_{i_1,j_1}\cdot h_{i,j}
& b_{i,j}\cdot b_{i_1,j_1} & =b_{i_1,j_1}\cdot b_{i,j}
\end{align*}}

\begin{proof}
It is a routine calculation in homology algebra that for the truncated polynomial algebra $\Gamma =T[x]$
with $|x|\equiv 0$ $mod$ $2$ and $x$ primitive,
\begin{align*}
 Ext_\Gamma(Z/p, Z/p) & =E[h]\otimes P[b]
\end{align*}
where $h\in Ext^1$ is represented in the cobar complex by $x$ and
$b\in Ext^2$ is represented by $\sum {p \choose m}/p(x^m \otimes x^{(p-m)})$
($b=h^2$ represented by $x\otimes x$ at the prime 2).
Notice that the $E_0$-term of the spectral sequence is isomorphic
to the cobar complex of $E^{*,M}(n,k)$. Then from (2.5)
we see that
$$H^{s,*,M}(E_0^{*,t,M}S(n,k), d_0)=Ext_{E^{*,*}(n,k)}^{s,t}(Z/p, Z/p)=
\bigotimes_{k\leqslant s}Ext_{T[t_s^{p^j}]}^{*,*}(Z/p, Z/p)$$
Thus the May's $E_1$-term
$$
E_1^{s,t,M}S(n,k)=E[h_{i,j}|k\leqslant i, j\in Z/n]\otimes
  P[b_{i,j}|k\leqslant i, j\in Z/n].
$$

    Notice that $d_0(t_i^{p^j}\cdot t_{i_1}^{p^{j_1}})=-t_i^{p^j}\otimes t_{i_1}^{p^{j_1}}
-t_{i_1}^{p^{j_1}}\otimes t_i^{p^j}$, we get $h_{i,j}h_{i_1,j_1}=-h_{i_1,j_1}h_{i,j}$.
In a similar way, one can proof that $h_{i,j}\cdot b_{i_1,j_1}  =b_{i_1,j_1}\cdot h_{i,j}$
and $b_{i,j}\cdot b_{i_1,j_1}  =b_{i_1,j_1}\cdot b_{i,j}$ (\cf [6] Lemma 3.4 and 3.8).
\end{proof}

\section{The first May differentials}

    Now suppose $k\leqslant n$, then $s_0\leqslant 2n$. From (2.2) and Lemma 2.4 one has
\begin{align}\tag\nr
d_1(h_{i,j}) = & -\sum_{k\leqslant r\leqslant i-k}h_{r,j}h_{i-r,j+r}
              & &\mbox{for $i\leqslant s_0$}\\
d_1(h_{i,j}) = & b_{i-n,j+n-1} & &\mbox{for $s_0<i$.} \notag
\end{align}
Thus for
$i>s_0-n$, $b_{i,j}$ is the boundary of the first May differentials.
Recall from [10] Theorem 4.3.22, in the cobar complex of  $BP/I_n$
one has
\begin{align}\tag\nr
d(b_{i,j}) = & \sum_{0<r<i}\left(b_{r,j}\otimes t_{i-r}^{p^{r+j+1}}-t_r^{p^{j+1}}\otimes b_{i-r,r+j}\right)
\end{align}
Thus for $i\leqslant s_0-n$, the first May differential $d_1(b_{i,j})=0$.
This implies:

{\bf Theorem \nr} {\it Let $k\leqslant n$ and $s_0$ be given in Lemma 2.4.
The $E_2-$term of the May spectral sequence is isomorphic to the
cohomology of
$$\w{E}_1^{*,*,*}S(n,k)=E[h_{i,j}|k\leqslant i\leqslant s_0, j\in Z/n]\otimes
         P[b_{i,j}|k\leqslant i\leqslant s_0-n, j\in Z/n].$$
The first May differential are given by
\begin{align*}
d_1(h_{i,j}) = & -\sum_{k\leqslant r\leqslant i-k}h_{r,j}h_{i-r,r+j}
         & &\mbox{for $i\leqslant s_0$} \\
d_1(b_{i,j}) = & 0 & & \mbox{for $k\leqslant i\leqslant s_0-n$}
\end{align*}
At the prime $p=2$, $s_0=2n$. The reduced May $E_1$-term becomes
$$\w{E}_1^{*,*,*}S(n,k)= E[h_{i,j}|n< i \leqslant 2n,\ j\in Z/n]
  \otimes P[h_{i,j}|k\leqslant i\leqslant n, j\in Z/n]$$
and the first May differential of $h_{2n,j}$ is given by
$$d_1(h_{2n,j})=-\sum_{k\leqslant i\leqslant 2n-k}h_{i,j}h_{2n-i,i+j} + h_{n,j+n-1}^2$$
}

\begin{proof}
We define a filtration in the May's $E_1$-term
$$E_1^{*,*,*}S(n,k)=E[h_{i,j}|k\leqslant i,\ j\in Z/n]\otimes P[b_{i,j}|k\leqslant i,\ j\in Z/n]$$
as follows: for each $s\geqslant k$, define
\begin{align*}
F^s(n,k) = &
\left\{ \begin{array}{lll}
          E[h_{i,j}|k\leqslant i\leqslant s] &
                  \mbox{for $k\leqslant s\leqslant n+k-1$}\\
          E[h_{i,j}|k\leqslant i\leqslant s]\otimes
          P[b_{i,j}|k\leqslant i\leqslant s-n] &
                  \mbox{for $n+k-1<s$.}
        \end{array}
\right.
\end{align*}
From (3.1), we see that for each $s\geqslant k$, $F^s(n,k)$ is a sub-complex
of $E_1^{*,*,*}S(n,k)$ that satisfies
$$F^{s_0}(n,k)=\w{E}_1^{*,*,*}(n,k)=E[h_{i,j}|k\leqslant i\leqslant s_0, j\in Z/n]\otimes
         P[b_{i,j}|k\leqslant i\leqslant s_0-n, j\in Z/n]$$
and
$$F^k(n,k)\hookrightarrow F^{k+1}(n,k) \hookrightarrow \cdots
  \hookrightarrow F^s(n,k) \hookrightarrow F^{s+1}(n,k) \hookrightarrow
  \cdots \hookrightarrow E_1^{*,*,*}(n,k).$$
Indeed, for $s>n+k-1$,
\[
F^s(n,k)=F^{s-1}(n,k)\bigotimes \left(E[h_{s,j}|j\in Z/n]\otimes P[b_{s-n,j}|j\in Z/n]\right).
\]

    For $s>s_0$ one has $d_1(h_{s,j})=b_{s-n,j+n-1}$. Thus
\[
E[h_{s,j}|j\in Z/n]\otimes P[b_{s-n,j}|j\in Z/n]
\]
is a sub-complex of $F^s(n,k)$
whose cohomology is $Z/p$ concentrated at dimensional 0.
This implies
$$H^*F^{s_0}(n,k)\cong H^*F^{s_0+1}(n,k) \cong \cdots  \cong H^*F^s(n,k) \cong
  \cdots \cong H^*E_1^{*,*,*}S(n,k)$$

  At prime $p=2$, $s_0=\left[\displaystyle{\frac{2\times2n}{2}}\right]=2n>n+k-1$. The first
May differentials are deduced from (2.2).
\end{proof}

    As a corollary one can easily see that if $\displaystyle{\frac{2pn+p-2}{2(p-1)}}\leqslant n+k-1$,
then the reduced May's $E_1$-term becomes
$$\w{E}_1^{*,*,*}S(n,k)=E[h_{i,j}|k\leqslant i\leqslant n+k-1, j\in Z/n].$$

{\bf Theorem \nr} {\it If $\displaystyle{\frac{2pn+p-2}{2(p-1)}}\leqslant n+k-1$, then the cohomology
of $S(n,k)$ is of dimensional $n^2$.}

\section{The higher May differentials in the MSS for $S(n,k)$}

    From (3.2) we see that the first non-trivial May differential of $b_{i,j}$ appears
at
\begin{align}\tag\nr
d_r(b_{i,j})=\left\{ \begin{array}{lll}
                    0 & &
                     \mbox{if $i<2k$.}\\
                    b_{i-k,j}h_{k,j+k+1}-h_{k,j+1}b_{i-k,j+k} & &
                     \mbox{if $i\geqslant 2k$.}
                   \end{array}
           \right.
\end{align}

    In [6] (2.10) and (2.11), a collapse theorem is given for the higher May differentials
in the exterior part $E[h_{i,j}|i>0,\ j\geqslant 0]$ of the MSS for the steenrod algebra $A$ at
odd primes. In this section, we will give a similar collapse theorem for the
higher May differentials of $E_r^{*,*,*}S(n,k)$.

    Let $p$ be an odd prime. We define a Hopf algebra $T(n,k)$ as
\begin{align}\tag\nr
T(n,k)= & P[\xi_i|k\leqslant i\leqslant n+k-1].
\end{align}
The inner degree of $\xi_i$ is defined to be $|\xi_i|=2(p-1)(1+p+\cdots +p^{i-1})$ and
the structure map $\Delta: T(n,k)\longrightarrow T(n,k)\otimes T(n,k)$ acts on $\xi_i$
by
$$\Delta(\xi_i)=\xi_i\otimes 1 + \sum_{k\leqslant r\leqslant i-k}\xi_r\otimes \xi_{i-r}^{p^r}
   + 1\otimes\xi_i.$$

  There is a Hopf algebra reduction
homomorphism $\Phi: T(n,k)\rightarrow S(n,k)$ which
send $\xi_i$ to $t_i$. The image of $\Phi$ is $P[t_i|k\leqslant i\leqslant n+k-1]/(t_i^{p^n}-t_i)$
and $Ker\ \Phi$ is the idea generated by $(\xi_i^{p^n}-\xi_i)$. Further more the homomorphism
$\Phi$ also induces homomorphism in cobar complexes and cohomologies
\[
\Phi: Ext_{T(n,k)}^{s,*}(Z/p, Z/p) \longrightarrow Ext_{S(n,k)}^{s,*}(Z/p, Z/p)
\]

  Similar to that of definition 2.3, we set May filtration on $T(n,k)$ as
$$M(\xi_i^{p^j})=2i-1$$
and let $F^{*,M}T(n,k)$ be the sub-module of $T(n,k)$ generated by the elements with
May filtration $\leqslant M$. Then $E^{*,M}T(n,k)=F^{*,M}T(n,k)/F^{*,M-1}T(n,k)$ becomes
a bigraded Hopf algebra with the structure of
$$E^{*,*}T(n,k)=\bigotimes T[\xi_i^{p^j}|k\leqslant i\leqslant n+k-1,\ j\geqslant 0]$$
and $\Delta(\xi_i^{p^j})=\xi_i^{p^j}\otimes 1+1\otimes \xi_i^{p^j}$.

    Consider the cobar construction $C^{s,*}T(n,k)$ of $T(n,k)$.
Similarly for the generator $[\beta_1|\beta_2|\cdots | \beta_s]$ of $C^{s,*}T(n,k)$
define its May filtration as
$$M([\beta_1|\beta_2|\cdots | \beta_s])=M(\beta_1)+M(\beta_2)+\cdots + M(\beta_s)$$
and let $FC^{*,*,M}T(n,k)$ denote the sun-complex  generated by elements with May filtration
$\leqslant M$. We get a spectral sequence $\{E_r^{s,*,M}T(n,k),\ d_r\}$ with $E_0$-term
$$E_0^{*,*,M}T(n,k)=FC^{*,*,M}T(n,k)/FC^{*,*,M-1}T(n,k)$$
which is isomorphic to the cobar complex of $E^{*,M}T(n,k)$.  The $E_1$-term of this
spectral sequence is isomorphic to
\begin{align}\tag\nr
E_1^{*,*,*}T(k,n)= & E[h'_{i,j}|k\leqslant i\leqslant n+k-1, j\geqslant 0]\otimes
              P[b'_{i,j}|k\leqslant i\leqslant n+k-1, j\geqslant 0].
\end{align}

   Noticed that the reduction map $\Phi: T(n,k)\rightarrow S(n,k)$  is May filtration preserving,
it induces a homomorphism of May spectral sequences
$$\Phi: E_r^{*,*,*}T(n,k) \longrightarrow E_r^{*,*,*}S(n,k).$$

{\bf Theorem \nr} {\it The reduction map $\Phi: T(n,k)\longrightarrow S(n,k)$ induces a homomorphism
between May spectral sequences $\Phi: E_1^{*,*,*}T(n,k) \longrightarrow E_1^{*,*,*}S(n,k)$
which sends $h'_{i,j}$ and $b'_{i,j}$ to $h_{i,j}$ and $b_{i,j}$ respectively. It sends infinite cocycles of
$E_r^{*,*,*}T(n,k)$ to that of $E_r^{*,*,*}S(n,k)$.}

    Similar to [6] (2.10) and (2.11) we give a collapse theorem in the MSS for $T(n,k)$.
To the generators $h'_{i,j},\ b'_{i,j}\in E_1^{*,*,*}T(n,k)$
define their index as
\[
SI(h'_{i,j})=SI(b'_{i,j})=i.
\]
Given a monomial
$g=x_1x_2\cdots x_m  \in E_1^{*,*,*}T(n,k)$
where each $x_i$ is of the generators $h'_{i,j}$ or $b'_{i,j}$,
define its sum of index as
\begin{align}\tag\nr
SI(g)=SI(x_1)+SI(x_2)+\cdots +SI(x_m).
\end{align}

For example the sum of index of $h'_{4,0}h'_{3,0}b'_{2,1}$ is $9$

  We use $s(x)$ to denote the homological dimension of $x$.
Noticed that the May filtration of $h'_{i,j},\ b'_{i,j}$ satisfies
\begin{align*}
M(h'_{i,j})= & 2i-1 =2SI(h'_{i,j})-1=2SI(h'_{i,j})-s(h'_{i,j})\\
M(b'_{i,j})= & p(2i-1)>2SI(b'_{i,j})-2=2SI(b'_{i,j})-s(b'_{i,j})
\end{align*}
we see that for the monomial $g=x_1x_2\cdots x_m\in E_1^{s, *,*}T(n,k)$
of homological dimension $s$, its May filtration satisfies
\begin{align}\notag
M(g)= & M(x_1)+M(x_2)+\cdots+M(x_m)\\
\geqslant & 2SI(x_1)-s(x_1)+2SI(x_2)-s(x_2)+\cdots +2SI(x_m)-s(x_m)\tag\nr\\
= & 2SI(g)-s \notag
\end{align}
and the equality holds if and only if $g$ is a monomial in
$E[h'_{i,j}|k\leqslant i\leqslant n+k-1, j\geqslant 0]$.

   Given an integer $t=2(p-1)(c_0+c_1p+\cdots +c_mp^m)$ with $0\leqslant c_i<p$,
we define its sum of degree as
\begin{align*}
Sd(t)=& c_0+c_1+\cdots+c_m
\end{align*}
and for an element $g\in E_1^{*,*,*}T(n,k)$, express its inner
degree $|g|$ as $|g|=2(p-1)(c_0+c_1p+\cdots +c_mp^m)$, where $0\leqslant c_i<p$ and define its
sum of degree to be
\begin{align}\tag\nr
Sd(g)=Sd(|g|)=& c_0+c_1+\cdots +c_m.
\end{align}
Then from
\begin{align*}
|h'_{i,j}|= & 2(p-1)(p^j+p^{j+1}+\cdots p^{i+j-1}) \\
|b'_{i,j}|= & 2(p-1)(p^{j+1}+p^{j+2}+\cdots +p^{i+j})
\end{align*}
we see that $SI(h'_{i,j})=Sd(h'_{i,j})$, $SI(b'_{i,j})=Sd(b'_{i,j})$.
But for the reason of the $p$-adic numbers one has
\begin{align}\tag\nr
SI(x_1x_2\cdots x_s)\geqslant & Sd(x_1x_2\cdots x_m).
\end{align}

{\bf Theorem \nr} {\it  In the May spectral sequence for $T(n,k)$,
\begin{enumerate}
\item If the inner degree $t=2(p-1)(c_0+c_1p+\cdots c_mp^m)$
and the May filtration
\[
M<2Sd(t)-s=2(c_0+c_1+\cdots +c_m)-s,
\]
then the May's $E_1$-term $E_1^{s,t,M}T(n,k)=0$.
\item If a cocycle $g\in E[h'_{i,j}|k\leqslant i\leqslant n+k-1, j\geqslant 0]$ in
the exterior part of May's $E_1$-term
satisfies $SI(g)=Sd(g)$, then it is an infinite cocycle in
the MSS for $T(n,k)$ and $\Phi(g)$ is an infinite cocycle in the MSS
for $S(n,k)$.
\end{enumerate}}

\begin{proof} (1) follows from (4.6) and (4.8).

Suppose $g\in E[h'_{i,j}|k\leqslant i\leqslant n+k-1, j\geqslant 0]$ is a cocycle in the exterior part of
May's $E_1$-term $E_1^{s,t,M}T(n,k)$ that satisfies $SI(g)=Sd(g)$. Then its May filtration $M=2SI(g)-s=2Sd(t)-s$.
Consider the higher May differentials
\[
d_r: E_r^{s,t,M}T(n,k)\rightarrow E_r^{s+1,t,M-r}T(n,k),
\]
we see that $M-r<2Sd(t)-(s+1)$ for $r>1$. Thus the target $E_1^{s+1,t,M-r}T(n,k)$ and then
$E_r^{s+1,t,M-r}T(n,k)$ is zero.
\end{proof}

{\bf Example} Let $p\geqslant 5$.  The $E_2$-term of the May spectral sequence for $H^{s,*}S(4,2)$
is isomorphic to the homology of
\[
E[h_{2,j}, h_{3,j}, h_{4,j}, h_{5,j}|j\in Z/4]
\]
with first May differentials
\begin{align*}
d_1(h_{2,j}) = & 0, & d_1(h_{3,j}) = & 0 \\
d_1(h_{4,j}) = & h_{2,j}h_{2,j+2}, & d_1(h_{5,j}) = & h_{2,j}h_{3,j+2} + h_{3,j}h_{2,j+3}.
\end{align*}
So $h_{5,0}h_{4,0}h_{3,0}h_{2,0}$ is a cohomology class in May's $E_2$-term.

  To prove that $h_{5,0}h_{4,0}h_{3,0}h_{2,0}$ is a infinite cocycle in the MSS
$E_r^{4,*}S(4,2)$, consider the MSS for $T(4,2)$.
$h_{5,0}h'_{4,0}h'_{3,0}h'_{2,0}$ is a
4-dimensional cocycle in the
exterior part of May's $E_1$-term $E_1^{4,t,M}T(4,2)$.
\[
deg(h'_{5,0}h'_{4,0}h'_{3,0}h'_{2,0})=2(p-1)(4+4p+3p^2+2p^3+p^4),
\]
\[
SI(h'_{5,0}h'_{4,0}h'_{3,0}h'_{2,0})=14=Sd(h'_{5,0}h'_{4,0}h'_{3,0}h'_{2,0}).
\]
Thus it is an infinite cocycle
in the MSS for $T(4,2)$ and $h_{5.0}h_{4,0}h_{3,0}h_{2,0}=\Phi(h'_{5,0}h'_{4,0}h'_{3,0}h'_{2,0})$ is an
infinite cocycle in the MSS for $S(4,2)$.

\section{The cohomology of $S(n,n)$ at $p=2$ and of $S(3,2)$ at $p=3$}

    As an application of Theorem 3.3, we will compute $H^{*,*}S(n,n)$ at $p=2$,
$H^{*,*}S(3,2)$ at prime $p=3$ and $H^{*,*}S(4,2)$ at prime $p\geqslant 5$ in this section.

\subsection{The cohomology of $S(n,n)$ at prime two}

    Consider the cohomology of $S(n,n)$ at $p=2$. The reduced May¡®s $E_1$-term
becomes
$$\w{E}_1^{*,*,*}S(n,n)=P[h_{n,j}|j\in Z/n]\otimes E[h_{s,j}|n<s\leqslant 2n, j\in Z/n]$$
(\cf Theorem 3.3). Noticed that the only non-trivial first May differential is
\begin{align}\tag\nr
d_1(h_{2n,j}) = & h_{n,j-1}^2+h_{n,j}^2,
\end{align}
We see that the $E_2$-term is the tensor product of
$E[h_{s,j}|n<s<2n]$ and the cohomology of
\[
\left\{P[h_{n,j}|j\in Z/n]\otimes E[h_{2n,j}|j\in Z/n],\ d_1\right\}.
\]

{\bf Lemma \nr} {\it The May's $E_2$-term $E_2^{*,*,*}S(n,n)$ at $p=2$ is isomorphic to the tensor product of
$E[h_{s,j}|n<s<2n, j\in Z/n]$ and $E[h_{n,j}, \rho_{2n}|j\in Z/n]\otimes P[h_{n,n-1}]$,
where $\rho_{2n}=\sum_{0\leqslant j<n} h_{2n,j}$ and $h_{n,j}^2=h_{n,n-1}^2$.}

\begin{proof}
We define $b_{n,j}=h_{n,j}^2+h_{n,j+1}^2$ for $0\leqslant j \leqslant n-2$ and define
$b_{n,n-1}=h_{n,n-1}^2$. It is easy to see that
$P[h_{n,j}|j\in Z/n]$ could be divided as the tensor product of
$P[b_{n,j}|0\leqslant j<n]$ and $E[h_{n,j}|j\in Z/n]$ as $Z/2$-modules.
\[
P[h_{n,j}|0\leqslant j\leqslant n-1]=
 P[b_{n,j}|0\leqslant j\leqslant n-1]\otimes E[h_{n,j}|0\leqslant j\leqslant n-1].
\]

   From (5.1) we see that
$$
d_1(h_{2n,j})=\left\{\begin{array}{lll}
                      b_{n,j-1} & & \mbox{if $1\leqslant j<n$} \\
                      \sum_{0\leqslant i<n-2}b_{n,i} & & \mbox{if $j=n$.}
                     \end{array}
              \right.
$$
The cohomology of
$\left\{P[h_{n,j}|j\in Z/n]\otimes E[h_{2n,j}|j\in Z/n],\ d_1\right\}$
is isomorphic to the tensor product of $E[h_{n,j}|j\in Z/n]$ and the cohomology of
$$P[b_{n,j}|0\leqslant j<n]\otimes E[h_{2n,j}|j\in Z/n].$$
The generator of $P[b_{n,j}|j\in Z/n]$ are denoted as
$$b_{n,i_1}^{s_1}b_{n,i_2}^{s_2}\cdots b_{n,i_m}^{s_m}$$
subject to $s_i>0$, $0\leqslant i_1<i_2<\cdots i_m<n$ and
the generators of $E[h_{2n,j}|j\in Z/n]$ are denoted as
$$h_{2n,j_k}\cdots h_{2n,j_2}h_{2n,j_1}$$
subject to $n\geqslant j_k>\cdots >j_2>j_1>0$.

    For the generators of $E[h_{2n,j}|0<j\leqslant n]\otimes P[b_{n,j}|0\leqslant j<n]$
described as above, one has
\begin{enumerate}
\item For $j_1>i_1+1$,
  $$h_{2n,j_k}\cdots h_{2n,j_2}h_{2n,j_1}b_{n,i_1}^{s_1}b_{n,i_2}^{s_2}\cdots b_{n,i_m}^{s_m}$$
is the leading term of the first May differential
$$d_1(h_{2n,j_k}\cdots h_{2n,j_2}h_{2n,j_1}h_{2n,i_1-1}
    b_{n,i_1}^{s_1-1}b_{n,i_2}^{s_2}\cdots b_{n,i_m}^{s_m}).$$

    While for $i_1<n-1$,
$$b_{n,i_1}^{s_1}b_{n,i_2}^{s_2}\cdots b_{n,i_m}^{s_m}$$
is the leading term of the first May differential
$$d_1(h_{2n,i_1+1}b_{n,i_1}^{s_1-1}b_{n,i_2}^{s_2}\cdots b_{n,i_m}^{s_m}).$$

\item For $j_1\leqslant i_1+1$ and $j_1<n$, the leading term of the first May differential
    $$d_1(h_{2n,j_k}\cdots h_{2n,j_2}h_{2n,j_1}b_{n,i_1}^{s_1}b_{n,i_2}^{s_2}\cdots b_{n,i_m}^{s_m})$$
is
  $$h_{2n,j_k}\cdots h_{2n,j_2}b_{2n,j_1-1}b_{n,i_1}^{s_1}b_{n,i_2}^{s_2}\cdots b_{n,i_m}^{s_m}.$$
\item For $j_1=n=i_1+1$,
    $$d_1(h_{2n,n}b_{n,n-1}^{s_1}\cdot x)=\sum_{i=0}^{n-2}b_{n,i}b_{n,n-1}^{s_1}\cdot x
           =d_1(\sum_{i=0}^{n-2}h_{2n,i+1}b_{n,n-1}^{s_1}\cdot x)$$
\end{enumerate}
    Thus the cohomology of $E[h_{2n,j}|j\in Z/n]\otimes P[b_{n,j}|j\in Z/n]$ is isomorphic to
$E[\rho_{2n}|j\in Z/n]\otimes P[b_{n,n-1}]$, where $\rho_{2n}=\sum_{0\leqslant i<n}h_{2n,i}$.
The Lemma follows.
\end{proof}

{\bf Theorem \nr} {\it The May $E_\infty$-term $E_\infty^{*,*,*}S(n,n)$ is isomorphic to
its $E_2$-term. Thus the cohomology of $S(n,n)$ at prime $2$ isomorphic to the tensor product of
$E[h_{s,j}|n<s<2n]$ and $E[h_{n,j}, \rho_{2n}|j\in Z/p]\otimes P[h_{n,0}]$}

\begin{proof}
It is easy to see from (2.2) that for $n\leqslant s<2n$, $h_{s,j}$ is an infinite cocycle.
From $d(t_{2n}+t_{2n}^2 +\cdots t_{2n}^{2^{n-1}})=0$ we get the infinite cocycle
$\rho_{2n}$. The Theorem follows.
\end{proof}

\subsection{The cohomology of $S(3,2)$ at prime 3}
    Now consider the cohomology of $S(3,2)$ at prime $p=3$. From Lemma 2.4 we see that the
$s_0=4$. Thus from Theorem 3.3 we see that the reduced May's $E_1$-term is
$$\w{E}_1^{*,*,*}S(3,2)=E[h_{2,j},\ h_{3,j},\ h_{4,j}|j\in Z/3]$$
and the first May differentials are given by
\begin{align}\notag
d_1(h_{2,j})= & 0 &
d_1(h_{3,j})= & 0 & & \mbox{and}\\
d_1(h_{4,j})= & -h_{2,j}h_{2,j+2}. \tag{\nr}
\end{align}
The May's $E_2$-term is isomorphic to
$$E_2^{*,*,*}S(3,2)=H^{*,*,*}\left(E[h_{2,j},\ h_{4,j}|j\in Z/3],\ d_1\right) \otimes
         E[h_{3,j}|j\in Z/3]$$

{\bf Lemma \nr} {\it  The May $E_2$-term for $H^*S(3,2)$ at $p=3$ has Poincare series
$(x^6+3x^5+6x^4+9x^3+6x^2+3x+1)(x+1)^3$. It is the tensor product of
$E[h_{3,j}|j\in Z/3]$ and the $Z/3$ module $\mathcal{C}$ generated by the following element.
$$\begin{array}{c}
\mbox{The generators of  $\mathcal{C}$}\\
\begin{tabular}{|l|l|l|l|l|l|l|l|}
\hline
\multicolumn{1}{|c|}{Dimension} &
\multicolumn{1}{|c|}{0} &
\multicolumn{1}{c|}{1} &
\multicolumn{1}{c|}{2} &
\multicolumn{1}{c|}{3} &
\multicolumn{1}{c|}{4} &
\multicolumn{1}{c|}{5} &
\multicolumn{1}{c|}{6} \\
\hline
Generators & $1$ & $h_{2,j}$ & $g_j$ & $l_j$  & $g_jg_{j+1}$ & $g_jl_{j+1}$ & $A$ \\
           &     &           & $k_j$ & $l'_j$ & $k_jk_{j+1}$ &              &     \\
           &     &           &       & $k_jh_{2,j}$ &        &              &     \\
\hline
\end{tabular}
\end{array}
$$
where $j\in Z/3$, $g_j=h_{4,j}h_{2,j}$, $k_j=h_{4,j}h_{2,j+2}$, $l_j=h_{4,j}h_{4,j+1}h_{2,j}$ and
\begin{align*}
l'_j= & h_{4,j}h_{4,j+1}h_{2,j+1}+h_{4,j+1}h_{4,j+2}h_{2,j} \\
A=    & h_{4,0}h_{4,1}h_{4,2}h_{2,0}h_{2,1}h_{2,2}=-g_0g_1g_2.
\end{align*}
}

\begin{proof}
    From (5.4), it is easy to see that $d_1(h_{4,j}h_{2,j})=0$, $d_1(h_{4,j}h_{2,j+2})=0$
and from $d_1(h_{4,j+1})=h_{2,j+1}h_{2,j+3}=h_{2,j+1}h_{2,j}$ we see that
$d_1(h_{4,j}h_{4,j+1}h_{2,j})=0$. These gives the cohomology classes
$g_j$, $k_j$ and $l_j$.  From
\begin{align*}
d_1(h_{4,j}h_{4,j+1}h_{2,j+1})= & -h_{2,j}h_{2,j+2}h_{4,j+1}h_{2,j+1}= -h_{4,j+1}h_{2,j}h_{2,j+1}h_{2,j+2}\\
d_1(h_{4,j+1}h_{4,j+2}h_{2,j})= & h_{4,j+1}h_{2,j+2}h_{2,j+4}h_{2,j} = h_{4,j+1}h_{2,j}h_{2,j+1}h_{2,j+2}
\end{align*}
we get $l'_j$. A routine computation shows that $H^*(E[h_{2,j},\ h_{4,j}|j\in Z/3])=\mathcal{C}$.
\end{proof}

{\bf Theorem \nr} {\it The May $E_2$-term for $H^{*,*}S(3,2)$ at $p=3$ is the $E_\infty$-term, thus
$H^{*,*}S(3,2)$ is the tensor product of $E[h_{3,j}|j\in Z/3]$ and $\mathcal{C}$.}

\begin{proof}
    It is easy to see that $h_{2,j}$ and $h_{3,j}$ are infinite cycles.
To prove that all the higher May differentials are trivial,
consider the May filtration of each generator in the $E_2$-term
and the differentials
\[
d_r: E_r^{s,t,M}S(3,2)\longrightarrow E_r^{s+1,t,M-r}S(3,2).
\]
One has
\begin{align*}
g_j \in & E_2^{2,*, 10}S(3,2) & k_j \in & E_2^{2,*, 10}S(3,2)\\
l_j \in & E_2^{3,*, 17}S(3,2) & l'_j \in & E_2^{3,*, 17}S(3,2)\\
h_{2,j} \in & E_2^{1,*, 3}S(3,2) & h_{3,j} \in & E_2^{1,0,5}S(3,2).
\end{align*}
The May filtration of $g_j$ and $k_j$ are 10. Beside, it is easy to check that each generator
in the 3rd dimension  $E_2^{3,*,*}S(3,2)$ listed as below
\[\begin{array}{ccccccccc}
l_j, & l'_j, & k_jh_{2,j}, & g_jh_{3,i}, & k_jh_{3,i}, & h_{3,i}h_{3,j}h_{2,k}, &
 h_{3,0}h_{3,1}h_{3,2}
\end{array}\]
has May filtration $\geqslant 10$. Thus $g_j$ and $k_j$
are infinite cycles. Similarly one can prove that $l_j$ and $l'_j$ are infinite cycles.
This complete the proof.
\end{proof}

\subsection{The cohomology of $S(4,2)$ at the primes $p>3$}
    In this case, $s_0=5$ and the reduced May's $E_1$-term is
$$\w{E}_1^{*,*,*}S(4,2)=E[h_{2,j}, h_{3,j}, h_{4,j}, h_{5,j}|j\in Z/4].$$

  To compute the $E_2$-term,
we set a filtration on the exterior algebra $E[h_{2,j}, h_{3,j}, h_{4,j}, h_{5,j}|j\in Z/4]$
as follows:
\[
F^k=\bigoplus_{0\leqslant r \leqslant k} Z/p\{h_{5,j_1}\cdots h_{5,j_r}\}
  \otimes E[h_{2,j}, h_{3,j}, h_{4,j}|j\in Z/4]
\]
where $h_{5,j_1}\cdots h_{5,j_r}$'s are the generators of the $r$-dimensional
module of the exterior algebra $E[h_{5,j}|j\in Z/4]$. This filtration gives raise
to a spectral sequence with
\[
E^k_0=F^k/F^{k-1}=Z/p\{h_{5,j_1}\cdots h_{5,j_k}\}\otimes E[h_{2,j}, h_{3,j}, h_{4,j}|j\in Z/4].
\]
The $E_1$-term of this spectral sequence is
\[
E_1^k=Z/p\{h_{5,j_1}\cdots h_{5,j_r}\} \otimes H^*E[h_{2,j}, h_{3,j}, h_{4,j}|j\in Z/4],
\]
and the differentials are given by
\[
d_r: E_r^k\longrightarrow E_r^{k-r}.
\]

  By a routine computation, we get

{\bf Theorem \nr} {\it
The cohomology of $E[h_{2,j}, h_{3,j}, h_{4,j}|j\in Z/4]$ is the
tensor product of $E[h_{3,j},\rho_0, \rho_1]$ and $\aleph$,
where
\begin{align*}
\rho_0= & h_{4,0} + h_{4,2}, & \rho_1= & h_{4,1}+h_{4,3},
\end{align*}
and $\aleph$ is
the direct sum of the modules
generated by the following cohomology classes:
$$\begin{array}{lllllllll}
1; & h_{2,j}; & e_j=h_{2,j}h_{2,j+1}, & g_j=h_{4,j}h_{2,j}; \\
 h_{2,j}g_{j+1}, & h_{2,j}g_{j+2}, & h_{2,j}g_{j+3}; \\
g_jg_{j+1}, & e_jg_{j+2}; & h_{2,j}g_{j+1}g_{j+2}; &
e_0g_2g_3
\end{array}$$
with $j\in Z/4$. Beside, we also have the following relations:
\begin{align*}
h_{2,i}h_{2,i+2}= & 0, & h_{2,i}g_{i+2}= & h_{2,i+2}g_i, & h_{2,i}g_{i+2}g_{i+3}= & h_{2,i+2}g_{i+3}g_i.
\end{align*}}

  With the add of a personal computer, we compute that

{\bf Theorem \nr} {\it
The cohomology of the exterior algebra $E[h_{2,j}, h_{3,j}, h_{4,j}, h_{5,j}]$ has Poincar\`{e} series
\[
(1+t)^4(1+6t+18t^2+59t^3+92t^4+176t^5+161t^6+176t^7+92t^8+59t^9+18t^{10}+6t^{11}+t^{12}).
\]
The ranks at each cohomological dimension are listed as
$$\begin{array}{cccccccccccccccccc}
0, & 1, & 2, & 3, & 4, & 5, & 6, & 7, & 8, & 9, &\cdots & 16,\\
1, & 10, & 48, & 171, & 461, & 976, & 1671, & 2303, &
  2558, & 2303, & \cdots & 1
\end{array}
$$
}

  From the collapse Theorem 4.9, we claim that the MSS for the cohomology of $S(4,2)$ collapse at $E_2$-term.

\end{document}